\documentclass[12pt]{article}
\usepackage{hyperref}
\usepackage{amsfonts}
\usepackage{enumerate}
\usepackage{amsthm}
\usepackage{amsmath}

\begin{document}

\begin{center}
{\large\bf \uppercase{The Samuelson's model as a singular discrete time system}}

\vskip.20in
Nicholas Apostolopoulos$^{1}$, Fernando Ortega$^{2}$ and\ Grigoris Kalogeropoulos$^{3}$\\[2mm]
{\footnotesize
$^{1}$National Technical University of Athens, Greece\\
$^{2}$ Universitat Autonoma de Barcelona, Spain\\
$^{3}$National and Kapodistrian University of Athens, Greece}
\end{center}

{\footnotesize
\noindent
\textbf{Abstract.} In this paper we revisit the famous classical Samuelson's multiplier-accelerator model for national economy. We reform this model into a singular discrete time system and study its solutions. The advantage of this study gives
a better understanding of the structure of the model and more deep and
elegant results.
\\[3pt]
{\bf Keywords}: Samuelson, macroeconomic, singular, system, difference equations
\\[3pt]

\vskip.2in

\section{Introduction}

Many authors have studied generalised discrete \& continuous time systems, see [1-19], and their applications especially in cases where the memory effect is needed including generalised discrete \& continuous time systems with delays, see [20-38]. Many of these results have already been extended to systems of differential \& difference equations with fractional operators, see [43-49]. 

Keynesian macroeconomics inspired the seminal work of Samuelson (1939), who actually introduced the business cycle theory. Although primitive and using only the demand point of view, the Samuelson's prospect still provides an excellent insight into the problem and justification of business cycles appearing in national economies. In the past decades, many more sophisticated models have been proposed by other researchers [20-38]. All these models use superior and more delicate mechanisms involving monetary aspects, inventory issues, business expectation, borrowing constraints, welfare gains and multi-country consumption correlations. 

Some of the previous articles also contribute to the discussion for the inadequacies of Samuelson's model. The basic shortcoming of the original model is: the incapability to produce a stable path for the national income when realistic values for the different parameters (multiplier and accelerator parameters) are entered into the system of equations. Of course, this statement contradicts with the empirical evidence which supports temporary or long-lasting business cycles. 

In this article, we propose an alternative view of the model by reforming it into a singular discrete time system. 

The paper is organized as follows. Section 2 provides a short review for the organization of the original model and in Section 3 we introduce the proposed reformulation into a system of difference equations. Section 4 investigates the solutions of the proposed system.      

\section{The original model}

The original version of Samuelson's multiplier-accelerator original model is based on the following assumptions:
\\\\
\textit{Assumption 2.1.} National income $T_k$ in year $k$, equals to the summation of three elements: consumption, $C_k$, private investment, $I_k$, and governmental expenditure $G_k$
\begin{equation}\label{eq11}
T_k=C_k+I_k+G_k.
\end{equation} 
\textit{Assumption 2.2.} Consumption $C_k$ in year $k$, depends on past income (only on last year's value) and on marginal tendency to consume, modeled with $a$, the multiplier parameter, where $0 < a < 1$,
\begin{equation}\label{eq12}
C_k=aT_{k-1}.
\end{equation}
\textit{Assumption 2.3.} Private investment $I_k$ in year $k$, depends on consumption changes and on the accelerator factor $b$, where $b>0$. Consequently,  $I_k$ depends on national income changes, 
\begin{equation}\label{eq13}
I_k=b(C_k-C_{k-1})=ab(T_{k-1}-T_{k-2}).
\end{equation}
\textit{Assumption 2.4.} Governmental expenditure  $G_k$ in year $k$, remains constant  
\[
G_k=\bar G.
\]
Hence, the national income is determined via the following second-order linear difference equation
\[
T_{k+2}-a(1+b)T_{k+1}+abT_k=\bar G.
\]
See [39-42] for the needed theory of difference equations that lead to the solution of the above equation.

\section{The reformulation - Singular Samuelson's model}

Let 
\[
Y_k=\left[\begin{array}{c}T_k\\C_k\\I_k\end{array}\right]
\]
Then \eqref{eq11} can be written as 
\[
0=-T_k+C_k+I_k+G_k,
\]
or, equivalently,
\[
\left[\begin{array}{ccc}0&0&0\end{array}\right]Y_{k+1}=\left[\begin{array}{ccc}-1&1&1\end{array}\right]Y_k+G_k.
\]
The equation \eqref{eq12} can be written as 
\[
C_{k+1}=aT_k
\]
or, equivalently,
\[
\left[\begin{array}{ccc}0&1&0\end{array}\right]Y_{k+1}=\left[\begin{array}{ccc}a&0&0\end{array}\right]Y_k.
\]
Finally \eqref{eq13} can be written as
\[
I_{k+1}=b(C_{k+1}-C_k).
\]
or, equivalently,
\[
-bC_{k+1}+I_{k+1}=-bC_k.
\]
or, equivalently,
\[
\left[\begin{array}{ccc}0&-b&1\end{array}\right]Y_{k+1}=\left[\begin{array}{ccc}0&-b&0\end{array}\right]Y_k.
\]
Hence the above expressions can be written in the following matrix form
\begin{equation}\label{eq1}
\begin{array}{cc}
FY_{k+1}=GY_k+V_k, & k= 2, 3,...,
\end{array}
\end{equation}
Where
\[
F=\left[\begin{array}{ccc}
0&0&0\\
0&1&0\\
0&-b&1\end{array}\right],\quad G=\left[\begin{array}{ccc}
-1&1&1\\
a&0&0\\
0&-b&0\end{array}\right],\quad V_k=\left[\begin{array}{c}
G_k\\0\\0\end{array}\right].
\]
Note that $F$ is singular (det$F$=0). Throughout the paper we will use in several parts matrix pencil theory to establish our results. A matrix pencil is a family of matrices $sF-G$, parametrized by a complex number $s$, see [46-53].
\\\\
\textbf{Definition 3.1.} Given $F,G\in \mathbb{R}^{r \times m}$ and an arbitrary $s\in\mathbb{C}$, the matrix pencil $sF-G$ is called:
\begin{enumerate}
\item Regular when  $r=m$ and  det$(sF-G)\neq 0$;
\item Singular when  $r\neq m$ or  $r=m$ and det$(sF-G)\equiv 0$.
\end{enumerate}
\textbf{Corollary 3.1.} The system \eqref{eq1} has always a \textsl{regular pencil} $\forall a,b$.\\\\
\textbf{Proof.} The determinant det$(sF-G)=s^2-a(b+1)s+ab\neq 0$. Hence from Definition 2.1, the pencil is regular. The proof is completed.
\\\\
The class of $sF-G$ is characterized by a uniquely defined element, known as the Weierstrass canonical form, see [50-57], specified by the complete set of invariants of $sF-G$. This is the set of elementary divisors of type  $(s-a_j)^{p_j}$, called \emph{finite elementary divisors}, where $a_j$ is a finite eigenvalue of algebraic multiplicity $p_j$ ($1\leq j \leq \nu$), and the set of elementary divisors of type $\hat{s}^q=\frac{1}{s^q}$, called \emph{infinite elementary divisors}, where $q$ is the algebraic multiplicity of the infinite eigenvalue. $\sum_{j =1}^\nu p_j  = p$ and $p+q=m$.
\\\\
From the regularity of $sF-G$, there exist non-singular matrices $P$, $Q$ $\in \mathbb{R}^{m \times m}$ such that 
\begin{equation}\label{eq3}
\begin{array}{c}PFQ=\left[\begin{array}{cc} I_p&0_{p,q}\\0_{q,p}&H_q\end{array}\right],
\\\\
PGQ=\left[\begin{array}{cc} J_p&0_{p,q}\\0_{q,p}&I_q\end{array}\right].\end{array}
\end{equation}
$J_p$, $H_q$ are appropriate matrices with $H_q$ a nilpotent matrix with index $q_*$, $J_p$ a Jordan matrix and $p+q=m$. With $0_{q,p}$ we denote the zero matrix of $q\times p$. The matrix $Q$ can be written as
\begin{equation}\label{eq4}
Q=\left[\begin{array}{cc}Q_p & Q_q\end{array}\right].
\end{equation}
$Q_p\in \mathbb{R}^{m \times p}$ and $Q_q\in \mathbb{R}^{m \times q}$. The matrix $P$ can be written as
\begin{equation}\label{eq5}
P=\left[\begin{array}{c}P_1 \\ P_2\end{array}\right].
\end{equation}
$P_1\in \mathbb{R}^{p \times r}$ and $P_2\in \mathbb{R}^{q \times r}$.
\\\\
The solution of system \eqref{eq1} is given by the following Theorem:
\\\\
\textbf{Theorem 3.1.}  (See [1-19]) We consider the system \eqref{eq1}. Since its pencil is always regular, its solution exists and for $k\geq 0$, is given by the formula
\[
    Y_k=Q_pJ_p^kC+QD_k.  
\]
Where $D_k=\left[
\begin{array}{c} \sum^{k-1}_{i=0}J_p^{k-i-1}P_1V_i\\-\sum^{q_{*}-1}_{i=0}H_q^iP_2V_{k+i}
\end{array}\right]$ and $C\in\mathbb{R}^p$ is a constant vector. The matrices $Q_p$, $Q_q$, $P_1$, $P_2$, $J_p$, $H_q$ are defined by \eqref{eq3}, \eqref{eq4}, \eqref{eq5}.

\section{Main Results}

In this section we will present our main results. We will provide the solution to the system \eqref{eq1} and consequently we will derive the sequence for the national income, the consumption and the private investment.
\\\\
\textbf{Theorem 4.1.} We consider the system \eqref{eq1}. Then in year $k$, National income $T_k$, Consumption $C_k$ and  private Investment $I_k$  are given by:
\[
\begin{array}{c} 
T_k=s_1^{k+1}c_1+s_2^{k+1}c_2+a\sum^{k-1}_{i=0}[(s_1^{k-1}+s_2^{k-1})]G_i,\\\\

C_k=a(s_1^kc_1+s_2^kc_2)+a^2\sum^{k-1}_{i=0}[(s_1^{k-i-1}+s_2^{k-i-1})]G_i,\\\\

I_k=s_1^{k}(s_1-a)c_1+s_2^{k}(s_2-a)c_2+a\sum^{k-1}_{i=0}[((s_1-a)s_1^{k-1}+(s_2-a)s_2^{k-1})]G_i
\end{array}
\]
\textbf{Proof.} From Corollary 3.1, the pencil $sF-G$ is always regular. Furthermore the pencil has one infinite eigenvalue and two finite:
\[
s_1=\frac{a(1+b)+\sqrt{a^2(1+b)^2-4ab}}{2},\quad s_2=\frac{a(1+b)-\sqrt{a^2(1+b)^2-4ab}}{2}.
\]
From Theorem 3.1, the solution of \eqref{eq1} is given by 
\[
    Y_k=Q_pJ_p^kC+Q\left[
\begin{array}{c} \sum^{k-1}_{i=0}J_p^{k-i-1}P_1V_i\\-\sum^{q_{*}-1}_{i=0}H_q^iP_2V_{k+i}
\end{array}\right].  
\]
Since we have one infinite eigenvalue we have 
\[
H_q=0
\]
and $J_p$ is the Jordan matrix of the two finite eigenvalues:
\[
    Y_k=Q_p
    \left[
\begin{array}{cc} s_1^k&0\\0&s_2^k
\end{array}\right]C+Q\left[
\begin{array}{c} \sum^{k-1}_{i=0}J_p^{k-i-1}P_1V_i\\0
\end{array}\right].  
\]
The matrix $Q_p$ has the two eigenvectors of the two finite eigenvalues:
\[
Q_p= \left[
\begin{array}{cc} s_1&s_2\\a&a\\s_1-a&s_2-a
\end{array}\right],
\]
 while $Q_q$ is the eigenvector of the infinite eigenvalue:
 \[
 Q_q= \left[
\begin{array}{c} 1\\0\\0
\end{array}\right].
 \]
 Hence:
 \[
 Q=\left[
\begin{array}{ccc} s_1&s_2&1\\a&a&0\\s_1-a&s_2-a&0
\end{array}\right]
 \]
 and the solution $Y_k$ takes the form:
\[
    Y_k=
    \left[
\begin{array}{cc} s_1&s_2\\a&a\\s_1-a&s_2-a
\end{array}\right]
    \left[
\begin{array}{cc} s_1^k&0\\0&s_2^k
\end{array}\right]C+
\]
\[
\left[
\begin{array}{ccc} s_1&s_2&1\\a&a&0\\s_1-a&s_2-a&0
\end{array}\right]
\left[
\begin{array}{c} \sum^{k-1}_{i=0} \left[
\begin{array}{cc} s_1^{k-i-1}&0\\0&s_2^{k-i-1}
\end{array}\right]P_1V_i\\0
\end{array}\right].  
\]
Finally, where $P_1$ is the matrix which contains the right eigenvectors of the finite eigenvalues
\[
  P_1=
    \left[
\begin{array}{ccc} a&1&\frac{a}{s_1}\\a&1&\frac{a}{s_2}
\end{array}\right].
\]
Hence
\[
    Y_k=
    \left[
\begin{array}{cc} s_1&s_2\\a&a\\s_1-a&s_2-a
\end{array}\right]
    \left[
\begin{array}{cc} s_1^k&0\\0&s_2^k
\end{array}\right]C+
\]
\[
\left[
\begin{array}{ccc} s_1&s_2&1\\a&a&0\\s_1-a&s_2-a&0
\end{array}\right]\left[
\begin{array}{c} \sum^{k-1}_{i=0}\left[
\begin{array}{cc} s_1^{k-i-1}&0\\0&s_2^{k-i-1}
\end{array}\right]\left[
\begin{array}{ccc} a&1&\frac{a}{s_1}\\a&1&\frac{a}{s_2}
\end{array}\right]\left[
\begin{array}{c} G_i\\0\\0
\end{array}\right]\\0
\end{array}\right],
\]
or, equivalently,
\[
    Y_k=
    \left[
\begin{array}{c} 
s_1^{k+1}c_1+s_2^{k+1}c_2+a\sum^{k-1}_{i=0}[(s_1^{k-1}+s_2^{k-1})]G_i\\

a(s_1^kc_1+s_2^kc_2)+a^2\sum^{k-1}_{i=0}[(s_1^{k-i-1}+s_2^{k-i-1})]G_i\\

s_1^{k}(s_1-a)c_1+s_2^{k}(s_2-a)c_2+a\sum^{k-1}_{i=0}[((s_1-a)s_1^{k-1}+(s_2-a)s_2^{k-1})]G_i

\end{array}\right],
\]
or, equivalently,
\[
   \left[\begin{array}{c}T_k\\C_k\\I_k\end{array}\right]=
    \left[
\begin{array}{c} 
s_1^{k+1}c_1+s_2^{k+1}c_2+a\sum^{k-1}_{i=0}[(s_1^{k-1}+s_2^{k-1})]G_i\\

a(s_1^kc_1+s_2^kc_2)+a^2\sum^{k-1}_{i=0}[(s_1^{k-i-1}+s_2^{k-i-1})]G_i\\

s_1^{k}(s_1-a)c_1+s_2^{k}(s_2-a)c_2+a\sum^{k-1}_{i=0}[((s_1-a)s_1^{k-1}+(s_2-a)s_2^{k-1})]G_i

\end{array}\right].
\]
The proof is completed.
\subsection*{Initial Conditions}
We assume system \eqref{eq1} and the known initial conditions (IC): $Y_{2}$.
\\\\
\textbf{Definition 4.1.} Consider the system \eqref{eq1} with known IC. Then the IC are called consistent if there exists a solution for the system \eqref{eq1} which satisfies the given conditions.
\\\\
\textbf{Proposition 4.2.} (See [1-19]) The IC of system \eqref{eq1} are consistent if and only if 
\[
Y_2\in colspanQ_p +QD_2.
\]
\textbf{Proposition 4.3.} (See [1-19]) Consider the system \eqref{eq1} with given IC. Then the solution for the initial value problem is unique if and only if the IC are consistent. Then, the unique solution is given by the formula
\[
    Y_k=Q_pJ_p^kZ^p_{2}+QD_k.  
\]
where $D_k=\left[
\begin{array}{c} \sum^{k-1}_{i=0}J_p^{k-i-1}P_1V_i\\-\sum^{q_{*}-1}_{i=0}H_q^iP_2V_{k+i}
\end{array}\right]$ and $Z^p_2$ is the unique solution of the algebraic system $Y_2=Q_pZ^p_2+D_2$.
\\\\
\textbf{Proposition 4.3.} The reformulation - Singular Samuelson's model has always a unique solution for given initial conditions\\\\
\textbf{Proof.} The reformulation - Singular Samuelson's model has always a unique solution for given initial conditions is a singular system given by \eqref{eq4}. For $k=2$ we get:
\[
  Y_2= \left[\begin{array}{c}T_2\\C_2\\I_2\end{array}\right],
  \]
  or, equivalently,
  \[
   Y_2 =
    \left[\begin{array}{c}T_2\\aT_1\\ab(T_1-T_0)\end{array}\right].
\]
or, equivalently,
  \[
   Y_2 =
    \left[\begin{array}{c}1\\0\\0\end{array}\right]T_2+\left[\begin{array}{c}0\\1\\b\end{array}\right]aT_1+\left[\begin{array}{c}0\\0\\-b\end{array}\right]aT_0.
\]
However
\[
colspanQ_p +QD_2=<\left[\begin{array}{c}0\\1\\b\end{array}\right],\left[\begin{array}{c}0\\0\\-b\end{array}\right]>+\left[\begin{array}{c}1\\0\\0\end{array}\right]
\]
and hence from Proposition 4.1, the IC of the reformulation - Singular Samuelson's model are always consistent and from Proposition 4.2, reformulation - Singular Samuelson's model has a unique solution for given IC. The proof is completed.


\begin{thebibliography}{00}
 


\bibitem{} Apostolopoulos, N., Ortega, F. and Kalogeropoulos, G., \emph{On stability of generalised systems of difference equation with non-consistent initial conditions}. arXiv preprint arXiv:1612.04120 (2016).

\bibitem{} Apostolopoulos, N., Ortega, F. and Kalogeropoulos, G., \emph{A boundary value problem of a generalised linear discrete time system with no solutions and infinitely many solutions.}. arXiv preprint arXiv:1610.08277 (2016).

\bibitem{} Apostolopoulos, N., Ortega, F. and Kalogeropoulos, G., \emph{The case of a generalised linear discrete time system with infinite many solutions.}. arXiv preprint arXiv:1610.00927 (2016).

\bibitem{}
I. Dassios, \emph{On a boundary value problem of a class of generalized linear discrete time systems,} Advances in Difference Equations, Springer, 2011:51 (2011).

\bibitem{}
I.K. Dassios, \emph{On non-homogeneous linear generalized linear discrete time systems}, Circuits systems and signal processing, Volume 31, Number  5, 1699-1712 (2012).
\bibitem{}
I. Dassios, \emph{On solutions and algebraic duality of generalized linear discrete time systems}, Discrete Mathematics and Applications, Volume 22, No. 5-6, 665--682 (2012).

\bibitem{} I. Dassios, \emph{On stability and state feedback stabilization of singular linear matrix difference equations,} Advances in difference equations, 2012:75 (2012).

\bibitem{} I. Dassios, \emph{On robust stability of autonomous singular linear matrix difference equations,} Applied Mathematics and Computation, Volume 218, Issue 12, 6912--6920 (2012).

\bibitem{}
I.K. Dassios, G. Kalogeropoulos, \emph{ On a non-homogeneous singular linear discrete time system with a singular matrix pencil }, Circuits systems and signal processing, Volume 32, Issue 4, 1615--1635 (2013).

\bibitem{}I. Dassios, G. Kalogeropoulos, \emph{On the relation between consistent and non consistent initial conditions of singular discrete time systems,} Dynamics of continuous, discrete and impulsive systems Series A: Mathematical Analysis, Volume 20, Number 4a, pp. 447--458 (2013).

\bibitem{}Dassios I., \emph{On a Boundary Value Problem of a Singular Discrete Time System with a Singular Pencil}, Dynamics of continuous. Discrete and Impulsive Systems Series A: Mathematical Analysis, 22(3): 211-231 (2015).

\bibitem{} I. K. Dassios, K. Szajowski, \emph{ Bayesian optimal control for a non-autonomous stochastic discrete time system,} Applied Mathematics and Computation , Volume 274, 556--564 (2016).

\bibitem{} I. K. Dassios, K. Szajowski, \emph{A non-autonomous stochastic discrete time system with uniform disturbances.} arXiv preprint arXiv:1612.05044, 2016.

\bibitem{}
I.K. Dassios, \emph{Homogeneous linear matrix difference equations of higher order: regular case}, Bull. Greek Math. Soc. 56, 57-64 (2009).

\bibitem{}Kalogeropoulos, Grigoris, and Charalambos Kontzalis. \emph{Solutions of Higher Order Homogeneous Linear Matrix Differential Equations: Singular Case.} arXiv preprint arXiv:1501.05667 (2015).


\bibitem{} F. L. Lewis; \emph{A survey of linear singular systems}, Circuits Syst. Signal Process. 5, 3-36, (1986).

\bibitem{} F.L. Lewis;\emph{ Recent work in singular systems}, Proc. Int. Symp. Singular systems, pp. 20-24, Atlanta, GA, (1987).

\bibitem{} F. L. Lewis;\emph{ A review of 2D implicit systems, Automatica (Journal of IFAC)}, v.28 n.2, p.345-354, (1992).

\bibitem{Ver}
L. Verde-Star; \emph{Operator identities and the solution of
linear matrix difference and differential equations}, Studies in
Applied Mathematics 91 (1994), pp. 153-177.


\bibitem{} Apostolopoulos, N., Ortega, F. and Kalogeropoulos, G., 2015. \emph{Causality of singular linear discrete time systems}. arXiv preprint arXiv:1512.04740.

\bibitem{Chari94}Chari, V. V., Optimal Fiscal Policy in a Business Cycle Model. {\em  Journal of Political Economy, Vol. 102, issue 4, p. 52-61}, (1994).
\bibitem{Chow85}Chow, G. C., A model of Chinese National Income Determination, {\em Journal of Political Economy, vol 93, No 4, p.782-792}, (1985).

\bibitem{} Dassios I, Kontzalis C: On the stability of equilibrium for a foreign trade model. Proceedings of the 32nd IASTED international conference 2012.

\bibitem{} Dassios I, Kontzalis C, Kalogeropoulos G: A stability result on a reformulated Samuelson economical model. Proceedings of the 32nd IASTED international conference 2012.

\bibitem{}I. Dassios, A. Zimbidis, \emph{The classical Samuelson's model in a multi-country context under a delayed framework with interaction,} Dynamics of continuous, discrete and impulsive systems Series B: Applications \& Algorithms, Volume 21, Number 4-5b pp. 261--274 (2014).

\bibitem{} I. Dassios, A. Zimbidis, C. Kontzalis. \emph{The Delay Effect in a Stochastic Multiplier-Accelerator Model.} Journal of Economic Structures 2014, 3:7.

\bibitem{}I. Dassios, G. Kalogeropoulos, \textit{On the stability of equilibrium for a reformulated foreign trade model of three
countries.} Journal of Industrial Engineering International, Springer, Volume 10, Issue 3, pp. 1-9 (2014). 10:71 DOI 10.1007/s40092-014-0071-9.

\bibitem{} I. Dassios, M. Devine. \emph{A macroeconomic mathematical model for the national income of a union of countries with interaction and tradel.} Journal of Economic Structures 2016, 5:18.



\bibitem{Day99} I. Dassios, A. Jivkov, A. Abu-Muharib and P. James. \emph{A mathematical model for plasticity and damage: A discrete calculus formulation}. Journal of Computational and Applied Mathematics, 2015. DOI 10.1016/j.cam.2015.08.017.
\bibitem{Dorf83}Dorf, R. C., Modern Control Systems. {\em  Addison-Wesley, 3rd Edition}, (1983).

\bibitem{} F. Milano; I. Dassios, Small-Signal Stability Analysis for Non-Index 1 Hessenberg Form Systems of Delay Differential-Algebraic Equations, Circuits and Systems I: Regular Papers, IEEE Transactions on 63(9):1521-1530 (2016).

\bibitem{} F. Milano; I. Dassios, Primal and Dual Generalized Eigenvalue Problems for Power Systems Small-Signal Stability Analysis, IEEE Transactions on Power Systems (2017).

\bibitem{Kuo87}Kuo, B. C., Automatic Control Systems. {\em  Prentice Hall, 5th Edition}, (1996).

\bibitem{Puu04}Puu, T., Gardini, L. and Sushko, I. A Hicksian
multiplier-accelerator model with floor determined by
capital stock. {\em Journal of Economic Behavior and
Organization, Vol. 56}, (2004).
\bibitem{Rosser00} Rosser, J. B., From Catastrophe to Chaos: A General
Theory of Economic Discontinuities. {\em Academic Publishers, Boston}, (2000).
\bibitem{Samuelson39} Samuelson, P. Interactions between the multiplier
analysis and the principle of acceleration. {\em  Review of
Economic Statistics}, (1939).
\bibitem{Westerhoff06}Westerhoff, F. H., Samuelson's multiplier-accelerator
model revisited. {\em Applied Economics Letters, Vol. 56, p. 86-92}, (2006).
\bibitem{Wincoop96} Wincoop, E. A multi-country real business cycle model. {\em  Scand. journal of economics vol. 23, p. 233--251}, (1996)

\bibitem{7}Ogata, K: Discrete Time Control Systems. Prentice Hall, (1987)

\bibitem{} W.J. Rugh; \emph{Linear system theory}, Prentice Hall International (Uk), London (1996).

\bibitem{} J.T. Sandefur; \emph{Discrete Dynamical Systems}, Academic Press, (1990).

\bibitem{}
A. P. Schinnar, \emph{The Leontief dynamic generalized inverse.} The Quarterly Journal of Economics 92.4 pp. 641-652 (1978).



\bibitem{}
I.K. Dassios, \emph{Optimal solutions for non-consistent singular linear systems of fractional nabla difference equations}, Circuits, Systems and Signal Processing, Springer, Volume 34, Issue 6, pp. 1769-1797 (2015). DOI 10.1007/s00034-014-9930-2

\bibitem{} I.K. Dassios, D. Baleanu, \emph{On a singular system of fractional nabla difference equations with boundary conditions,} Boundary Value Problems, 2013:148 (2013).

\bibitem {}
I.K. Dassios, D.I. Baleanu. \emph{Duality of singular linear systems of fractional nabla difference equations.} Applied Mathematical Modeling, Elsevier, Volume 39, Issue 14, pp. 4180-4195 (2015). DOI 10.1016/j.apm.2014.12.039

\bibitem{} I. Dassios, D. Baleanu, G. Kalogeropoulos, \emph{On non-homogeneous singular systems of fractional nabla difference equations,} Applied Mathematics and Computation, Volume 227, 112--131 (2014).

\bibitem{} I. Dassios, \emph{Geometric relation between two different types of initial conditions of singular systems of fractional nabla difference equations}, \textit{Math. Meth. Appl. Sci.}, 2015, doi: 10.1002/mma.3771.

\bibitem{} I. Dassios, \emph{Stability and robustness of singular systems of fractional nabla difference equations}. Circuits, Systems and Signal Processing (2016). doi:10.1007/s00034-016-0291-x

\bibitem{} I. Podlubny, \emph{Fractional Differential Equations, Mathematics in Science and Engineering},p. xxiv+340. Academic Press, San Diego, Calif, USA (1999).


\bibitem{ChYa}
H.-W. Cheng and S. S.-T. Yau; \emph{More explicit formulas for the matrix exponential}, Linear Algebra Appl. 262 (1997), pp. 131-163.

\bibitem{Ga} B.N. Datta; \emph{Numerical Linear Algebra and Applications}, Cole Publishing Company, 1995.

\bibitem{}
L. Dai, \emph{Singular Control Systems}, Lecture Notes in Control and information Sciences Edited by M.Thoma and A.Wyner (1988).

\bibitem{}
R. F. Gantmacher; \emph{The theory of matrices I, II}, Chelsea,
New York, (1959).

\bibitem{}
G. I. Kalogeropoulos; \emph{Matrix pencils and linear systems}, Ph.D Thesis, City University, London, (1985).
\bibitem{}Kontzalis, Charalambos P., and Panayiotis Vlamos. \emph{Solutions of Generalized Linear Matrix Differential Equations which Satisfy Boundary Conditions at Two Points.} Applied Mathematical Sciences 9.10 (2015): 493-505.

\bibitem{Leo} I. E. Leonard; \emph{The matrix exponential}, SIAM Review Vol. 38,
No. 3 (1996), pp. 507-512.

\bibitem{} G.W. Steward and J.G. Sun; \emph{Matrix Perturbation Theory}, Oxford University Press, (1990).


\end{thebibliography}
\end{document}